\documentclass[11pt,a4paper,openany]{article}

\usepackage{amsmath,amsfonts,ams}
\usepackage{latexsym}
\usepackage{amssymb}

\numberwithin{equation}{section}

\begin{document}
\title{Extensions of Hardy inequality }
\author{{Jia Yuan and Junyong Zhang
$^\dag$ } \\
        {\small The Graduate School of China Academy of Engineering Physics  }\\
        {\small P. O. Box 2101,\ Beijing,\ China,\ 100088 } \\
        {\small   jiajia2377@sohu.com$^\dag$}\\
        {\small  zhangjunyong111@sohu.com$^\dag$}
         \date{}
        }
\maketitle
\begin{abstract} In
this paper we prove sharp Hardy inequalities by using  Maximal
function theory. Our results improve and extend the well-known
results of G.Hardy \cite{Ha04}, T.Cazenave \cite {Ca03},
J.-Y.Chemin\cite {Ch06} and T.Tao\cite {TT06}.


\end{abstract}
\section{Introduction}

\quad\quad The initial Hardy-type inequality is the following
type, given by G.Hardy \cite{Ha04}.\\

 {\bf Proposition 1}\quad If $p>1, f(x)\geq0,$ and $F(x)=\int_0^x
 f(t)dt,$ then
 \begin{eqnarray}\label{1.1}
\int_0^{\infty}\bigg(\frac{F}{x}\bigg)^p dx<
\bigg(\frac{p}{p-1}\bigg)^p\int_0^{\infty}f^p dx ,
 \end{eqnarray}
 unless $f\equiv0.$ The constant is the best possible.\\

 \quad In general, for $1<p<n,$  we have
 \begin{eqnarray}\label{a1.1}\bigg\|\frac{f}{|x|}\bigg\|_p\leq
 \frac{p}{n-p}\|\nabla
 f\|_p.\end{eqnarray}
 It is an immediate result of the
following proposition when $p=q<n$. The proof of which is given by T.Cazenave
\cite {Ca03}.\\

 {\bf Proposition 2}\quad Let $1\leq p<\infty.$ If
$q<n$ is such that $0\leq q \leq p$, then
$\frac{|u(\cdot)|^p}{|\cdot|^q}\in L^1(\R^n)$ for every $u\in
W^{1,p}(\R^n)$. Furthermore,
 \begin{eqnarray}\label{1.2}
\int_{\R^n}\frac{|u(\cdot)|^p}{|\cdot|^q}dx\leq\Big(\frac{p}{n-q}\Big)^q\|u\|_{L^p}^{p-q}\|\nabla
u\|_{L^p}^q,
 \end{eqnarray}
for every $u\in W^{1,p}(\R^n)$.

Our goal is to prove the following theorem:\\

{\bf Theorem}\quad Let $p>1, 0\leq s<{\frac n p}$, then there exists
a constant $C$  such that for  ${\forall u}\in \dot W^{s,p}(\R^n)$,
\begin{eqnarray}\label{1.3}
\int_{\R^n}\frac{|u(x)|^p}{|x|^{sp}}dx\leq C \|u\|_{\dot W^{s,p}}^p.
 \end{eqnarray}

{\bf Remark:}\quad $(i)$ \quad If $s=0$, it is obvious that the
result holds true for $\forall\
 {1<p<\infty}$.  Compare  (\ref{1.3}) with  (\ref{a1.1}), condition $sp<n$
 is essential, and  so (\ref{1.3})  can admit
 $p\ge n$.

\vskip0.2cm
$(ii)$ \quad  For  $\forall$  $0\leq s \leq 1$,  by interpolation,
we know that
$$ \|u\|_{\dot W^{s,p}}^p\leq C\|u\|_{L^p}^{p(1-s)}\|\nabla
u\|_{L^p}^{ps}. $$
Hence  (\ref{1.3}) implies  (\ref{1.2}) without considering the constant
if  we choose $s=\frac q p$. \\

\vskip0.2cm
 $(iii)$\quad
Without considering the constant $C$, the inequality is sharp. In
details, for $p=\frac{n}{s}$, if we take $f(x)\in C_c^{\infty}(\R^n)$
satisfying $f(x)=1$ for $|x|\leq1$ $f(x)=0$ for $|x|\geq2$, then
 the above theorem can not hold. In fact
$$\int_{\R^n}\frac{|f(x)|^p}{|x|^{n}}dx\geq{\int_{|x|\leq1}\frac{1}{|x|^{n}}dx}=\infty$$
while\\
$$\|f\|_{\dot W^{s,p}}^p<\infty.$$

\vskip0.2cm
$(iv)$ \quad For the case $p=2,  0\leq s<\frac{n}{p},$  the corresponding
result is\\
$$\int_{\R^n}\frac{|u(\cdot)|^2}{|\cdot|^{2s}}dx\leq C\|u\|_{\dot{H^s}}^2 .$$
Chemin proves it by duality method and Littlewood-Paley
decomposition, Tao uses frequency splitting technology giving
another proof. The details can be found in in J.-Y. Chemin\cite
{Ch06} and T.Tao\cite {TT06}. The two methods can't be used in our
case.  The method involved in our paper is different from those
given by Cazenave , Chemin and T.Tao; it relies on the method of
Riesz potential for homogeneous Sobolev spaces and Maximal
function theory  cited in Section 2.\\

 The theorem will be  proved in Section 3.  First we give some notations
and preliminary work.
\section{Preliminaries}

\quad First we introduce a kind of definition to homogeneous Sobolev
space by using Reisz potential, we can find the detail in
M.Stein\cite{ST70} and C. Miao
\cite{Mi04} .\\

{\bf { Definition 1}}\quad  Define Reisz potential $I_{\alpha} f$ by

$$I_{\alpha} f=(-\triangle)^{-\frac{\alpha}2}
f(x)=C_{n,\alpha}\int_{\R^n}|x-y|^{-n+{\alpha}}f(y)dy,$$ where
$C_{n,\alpha}$ is the constant dependent in $\alpha$ and $n$.
The norm of homogenous sobolev space is given by\\
$$\|f\|_{{\dot
W}^{s,p}}=\|(-\triangle)^{\frac{s}2} f(x)\|_p=\|I_{-s} f\|_p.$$

 {\bf{ Definition 2}}\quad
 Let $f\in \ L_{loc}({\R^n}),x\in\R^n$,
$B\subset\R^n$ be a sphere and $x\in B,$ we define
$$ Mf(x)=\sup_{x\in B}{\frac 1{|B|}}\int_B |f(y)|dy,\quad |B|=m(B) $$
$Mf$ is  called H-L maximal function, $M$ is called maximal
operator.\\

If $B(r,x)\subset\R^n$ is a sphere with it's center and radius $x,$
$r$, define
$$ Mf(x)=\sup_{r>0}{\frac 1{|B(r,x)|}}\int_{B(r,x)} |f(y)|dy$$
 $Mf$ is called H-L centered maximal function,\quad M is called centered
maximal operator.\\

As we  know, the two definitions are equivalent, maximal operator has the following properties:\\

{\bf{Lemma}}\quad Let $f(x)$ be a measurable function on $\R^n,$
then \\

 \quad (i)\quad If $f\in \ L^p({\R^n}),1\leq p\leq \infty,$ then for
 a.e.$x\in\R^n$, $Mf(x)< \infty.$\\

 \quad (ii)\quad Operator $M$ is weak (1,1) type , in detail,  $\forall
 {f\in\ L^1(\R^n)}, \alpha > 0$, we have\\
 $$ m\{{x|x\in\R^n;Mf(x)>\alpha}\}\leq \frac
 A\alpha\int_{\R^n}|F(x)|dx,$$
 where constant {A} depends on {n}.\\

 \quad (iii)\quad  Let $1<p\leq\infty$, then  operator $M$ is strong ($p$,$p$)
 type, that is $\forall
 {f\in\ L^p(\R^n)},$ we have ${Mf\in  L^p(\R^n)},$ with\\
 $$\|Mf\|_p\leq A_p\|f\|_p ,$$
 where $A_p$ is a constant only dependent on $n$
and $p$.\\

As an immediate result, we have the following corollary:\\

{\bf{Corollary}}\quad For $\forall
 {f\in\ L^{p'}(\R^n)},\forall {{p'}>q}$,we have
 $$ \|(M(|f|^q))^\frac 1 q\|_{p'}\leq C \|f\|_{p'}.$$
 In fact, as ${p'}>q$, that is $\frac{p'}{q}>1,$ by $\textbf{Lemma
},$ we have\\
$$\|(M(|f|^q))^\frac 1 q\|_{L_{p'}}=\|M(|f|^q)\|_{L^{\frac {p'} q
}}^\frac{1}{q}\leq C \||f|^q\|_{L^{\frac {p'} q}}^{\frac{1}{q}} = C
\|f\|_{L^{p'}.}  \quad \quad \Box $$

\section{Proof of the theorem}

Proof:\quad Let ${I_{-s} u}=f$, then $u=I_s f$.\;  Using the definition of
Sobolev space, it is sufficient to show that:\\
\begin{eqnarray}\label{3.1}
\bigg\|\frac{I_s f}{|x|^s}\bigg\|_p\leq C\|f \|_p.
\end{eqnarray}
Let
\begin{align*}
A f&=\frac{I_s f}{|x|^s}=\int_{\R_n} \frac{f(y)}{|x-y|^{n-s}|x|^s}dy \\
&=\int_{|x-y|\leq 100|x|}\frac{f(y)}{|x-y|^{n-s}|x|^s}dy +
\int_{|x-y|\geq100|x|}\frac{f(y)}{|x-y|^{n-s}|x|^s}dy \\
&={A_1 f + A_2 f}
\end{align*}
To prove (\ref{3.1}), we need only prove that both $A_1 $ and
$A_2 $ are strong $(p,p)$ type.\\
We consider $A_1 f$ first. Notice that $s > 0,$ we have
\begin{align*}
A_1 f&=\sum_{j\leq
0}\int_{{2^{j-1}100|x|}\leq|x-y|\leq{2^j}100|x|}\frac{|f(y)|}{|x-y|^{n-s}|x|^s}dy\\
&=\sum_{j\leq
0}\int_{|x-y|\sim{2^j}100|x|}\frac{|f(y)|}{|x-y|^{n-s}|x|^s}dy\\
&\leq \sum_{j\leq
0}\int_{|x-y|\sim{2^j}100|x|}\frac{|f(y)|}{{2^j}100|x|^{n-s}|x|^s}dy\\
&\leq\sum_{j\leq
0}\int_{|x-y|\leq{2^j}100|x|}\frac{|f(y)|}{2^{nj}|x|^s}dy\cdot
2^{js}\\
&\leq C{\sum_{j\leq 0}2^{js}}Mf\\
&\leq C' Mf
\end{align*}
as $p>1,$ from $\textbf{Lemma}$, we have $\|A_1 f\|_p\leq
C\|f\|_p.$\\

Now, we are in position to consider  $A_2 f $ . It is easy to see
$$A_2 f(x)=\int_{|x-y|\geq100|x|}\frac{f(y)}{|x-y|^{n-s}|x|^s}dy
\leq \int_{|y|\geq{99|x|}}\frac{f(y)}{|y|^{n-s}|x|^s}dy \triangleq B_2 f(x).$$

For $\forall {g\in\ L^{p'}(\R^n)}$,
\begin{align*}
\big(B_2
f(x),g(x)\big)&=\bigg(\int_{|y|\geq{99|x|}}\frac{f(y)}{|y|^{n-s}|x|^s}dy,g(x)\bigg)\\
&=\int_{\R^n}{\int_{|y|\geq{99|x|}}\frac{f(y)}{|y|^{n-s}|x|^s}dy
\cdot
g(x)dx}\\
&=\int_{\R^n}{\frac{1}{|y|^{n-s}}\int_{|x|\leq\frac{|y|}{99}}\frac{g(x)}{|x|^{s}}dx\cdot
f(y)dy}\\
&=\big(Tg(y),f(y)\big)
\end{align*}
where
$$Tg(y)={\frac{1}{|y|^{n-s}}\int_{|x|\leq\frac{|y|}{99}}\frac{g(x)}{|x|^{s}}dx.
}$$

To prove $B_2$ to be strong $(p, p)$ type, it is sufficient to show that  \\
\begin{eqnarray}\label{3.2}
\|Tg(y)\|_{p'}\leq C\|g\|_{p'}.
\end{eqnarray}

In fact, consider
\begin{align*}
\big(B_2f(x),g(x)\big)&=\big(Tg(y),f(y)\big)\\
&=\|Tg(y)\|_{p'}\|f(y)\|_p\\
&\leq C\|g\|_{p'}\|f(y)\|_p,
\end{align*}
we  get that $B_2$ is strong $(p,p)$ type
by definition of norm.  Hence, $A_2$ is also a strong $(p,p)$.
 Now let us
prove formula (\ref{3.2}).

 For $s>0,$ when ${sq'}<n$, that is
$q>\frac{n}{n-s},$ we have, by H\"{o}lder inequality
\begin{align*}
{|Tg(y)|} &\leq
{\frac{1}{|y|^{n-s}}\bigg(\int_{|x|\leq\frac{|y|}{99}}|g(x)|^qdx\bigg)^{\frac{1}{q}}
\bigg(\int_{|x|\leq\frac{|y|}{99}}\frac{1}{|x|^{sq'}}dx\bigg)^{\frac{1}{q'}}} \\
&\leq C\frac{1}{|y|^{n-s}}\cdot\bigg(\int_{|x|\leq\frac{|y|}{99}}|g(x)|^qdx\bigg)^{\frac{1}{q}}\cdot{|y|^{{(n-{sq'})}\frac{1}{q'}}}\\
&\leq C\frac{1}{|y|^{\frac{n}{q}}}\cdot\bigg(\int_{|x|\leq\frac{|y|}{99}}|g(x)|^qdx\bigg)^{\frac{1}{q}}\\
&\leq C(M(|g|^q))^{\frac{1}{q}}(x)
\end{align*}
for $\forall {p'\geq q>{\frac{n}{n-s}},}$  that is
$1<p<{\frac{n}{s},}$ followed by {\bf{Corollary}}, we get $T$ strong
($p'$,$p'$)
type. This proves the theorem.\quad \quad $\Box$ \\

 {\bf Remark:}\quad For $p=1$, using our method, we can only get
$A_1$ is
weak $(1,1)$ type.\\

{\bf Acknowledge:}\quad The author is grateful to Prof. Xiaochun Li
 and Prof. Changxing Miao for their valuable suggestions.
\begin{center}

\end{center}
\end{document}